\documentclass{amsart}
\usepackage{latexsym}
\newtheorem{theorem}{Theorem}
\newtheorem{dichotomy}{Dichotomy}
\newtheorem{Lemma}{Lemma}

\usepackage[all]{xy}
\usepackage{epsf}

\numberwithin{equation}{section}

\newcommand{\be}{\begin{equation}}
\newcommand{\ee}{\end{equation}}
\newcommand{\bes}{\begin{equation*}}
\newcommand{\ees}{\end{equation*}}
\newcommand{\e}{\epsilon}
\newcommand{\te}{\theta}
\newcommand{\T}{\Theta}

\begin{document}

\title{Characterizing Jacobians via flexes of the Kummer Variety}

\author{Enrico Arbarello, Igor Krichever, Giambattista Marini}
\date{January 2005}


\maketitle

\begin{abstract} Given an abelian variety $X$ and a point $a\in X$
we denote by $<a>$ the closure of the subgroup of $X$ generated by $a$.
Let $N=2^{g}-1$. We denote by $\kappa: X\to \kappa(X)\subset\mathbb P^N$ the map from $X$ to its Kummer variety.
We prove that an indecomposable abelian variety $X$
is the Jacobian of a curve if and only if there exists a point $a=2b\in X\setminus\{0\}$
such that $<a>$ is irreducible and $\kappa(b)$
is a flex of $\kappa(X)$.

\end{abstract}

\section*{Introduction}\label{s: Introduction}

Let us begin by briefly recalling a few aspects of the KP equation:

\be
\label{kp1}
3u_{yy}=\frac{\partial}{\partial
x}\left[4u_t-6uu_{x}-u_{xxx}\right]\,,\quad u=u(x,y,t)\,.
\ee
It admits the, so-called, zero-curvature representation (\cite{zakh,dr})
\be\label{zcurv}
\left[L-\frac{\partial}{\partial y}\,,\,\,
A-\frac{\partial}{\partial t}\right]=0\,.
\ee
\noindent
where $L$ and $A$ are differential operators of the form
\be
\label{LA}
\aligned
L=&\frac{\partial^2}{\partial x^2}+u\,,\\
A=&\frac{\partial^3}{\partial x^3}+\frac{3}{2}u\frac{\partial}{\partial x}+w\,,
\qquad w=w(x,y,t)\,.
\endaligned
\ee
Equation (\ref{zcurv}) is the compatibility condition for an over-determined
system
of the linear equations:
\be
\label{eigen}
\aligned
\left(L-\frac{\partial}{\partial y}\right)\psi=0\,,\\
\left(A-\frac{\partial}{\partial t}\right)\psi=0\,.
\endaligned
\ee
A solution $\psi =\psi(x,y,t;\epsilon)$ of equations (\ref{eigen}) having the
form

\be
\label{wave}
\psi(x,y,t;\epsilon)=e^{\frac{x}{\epsilon}+\frac{y}{\epsilon^2}+\frac{t}{\epsilon^3}}
\left(1+\xi_1\epsilon+\xi_2\epsilon^2+\dots\right)\,.
\ee
is called the \it wave function\rm.
Here $\epsilon$ is a formal parameter and $\xi_i=\xi_i(x,y,t)$.

In \cite{kri1}, \cite{kri2} the general algebraic-geometrical
construction of quasi-periodic solutions of two-dimensional soliton equations of
the
KP type was proposed. This construction is based on the
concept of the Baker-Akhiezer function $\psi(x,y,t,Q)$, which is
uniquely determined by its analytical properties  on an auxiliary
Riemann surface $C$ and a point $Q\in C$.  The corresponding analytical
properties  generalize the analytical properties of the Bloch
functions of ordinary finite-gap linear periodic Sturm-Liouville  operators
established by Novikov, Dubrovin, Matveev and Its (see \cite{NDM},
\cite{sol} and references therein; see also \cite{lax}, \cite{mckean}).

Let $C$ be an algebraic curve (smooth and connected) of positive genus $g$. Let $\varphi: C\to J(C)$ be the Abel-Jacobi map with base point $p_0\in C$. In terms of a local parameter $\epsilon$ around $p_0$ and vanishing at $p_0$,
a local lifting to $\mathbb C^g$ of the map $\frac{1}{2}\varphi$  can be written as
\be
\label{zeta}
\epsilon\mapsto\zeta(\epsilon)=U\e+V\e^2+W\e^3+\dots\,,\qquad U,V,W\in\mathbb C^g\,.
\ee
Let $(z_1,\dots z_g)$ be coordinates in $\mathbb C^g$ and set $z=(z_1,\dots z_g)$. Let $\tau$ be the 
normalized period matrix of $C$ and consider the corresponding Riemann theta function
\be
\theta(z,\tau)=\sum_{n\in\mathbb Z^g}\exp 2\pi i\left(\frac{1}{2}nz{}^tn+z{}^tn\right)\,.
\ee
Then 
\be
\label{utheta}
u(x,y,t)=2\frac{\partial^2}{\partial x^2}\log\te(xU+yV+tW+z)+c
\ee
is a solution of the KP equation (\ref{kp1}), for any $z\in \mathbb C^g$ and $c\in \mathbb C$.
In this setting, the wave function becomes the Baker-Akhiezer function
\be
\label{bak}
\psi(x,y,t;\epsilon,z)=e^\Lambda\cdot
\frac{\te(xU+yV+tW+\zeta(\e)+z)}{\te(xU+yV+tW+z)}\,,
\ee
where
\be
\Lambda={\frac{x}{\epsilon}+\frac{y}{\epsilon^2}+\frac{t}{\epsilon^3}+\e \Lambda_1+\e^2\Lambda_2}+\dots
\ee
with $\Lambda_1, \Lambda_2,\dots$ are linear forms in $x,y$ and $ t$ having as coefficients holomorphic functions in $z$. 
Writing $U=(U_1,\dots, U_g)$, and similarly $V$ and $W$, we introduce the vector fields
\be
\label{fields}
D_1=\sum U_i\frac{\partial}{\partial z_i}\,,\quad D_2=\sum V_i\frac{\partial}{\partial z_i}\,,\quad D_3=\sum W_i\frac{\partial}{\partial z_i}\,.
\ee
We now plug in the KP equation (\ref{kp1}),  the expression for  $u$  given in (\ref{utheta}). We get the equation
\be
\label{hirota}
\aligned
&D_1^4\te\cdot\te-4D_1^3\te\cdot D_1\te+3D_1^2\te\cdot D _1^2\te+
3D_2^2\te\cdot \te\\&-3D_2\te \cdot D_2\te
-3D_1D_3\te \cdot \te+3D_3\te\cdot D_1\te-d\te\cdot\te=0\,.
\endaligned
\ee
where $\te=\te(z)$ and $d\in \mathbb C$. This is the KP equation  in Hirota bilinear form. Now start from a general
principally polarized abelian variety $(X,\T)$ where, as usual, $\T=\{x\in X\,|\,\te(x)=0\}$. Given constant vector fields $D_1$ and $D_2$ as in (\ref{fields}), we may consider the subschemes
\be\aligned
D_1\T&=\{x\in\T\,\,|\,\,D_1\te(x)=0\}\subset\T\,,\\
(D_1^2\pm D_2)\T&=\{x\in D_1\T\,\,|\,\,(D_1^2\pm D_2)\te(x)=0\}\subset D_1\T\,.
\endaligned
\ee
Clearly the KP equation (\ref{hirota}) implies a Weil-type relation
\be
\label{weil}D_1\T\,\,\subset \,\,(D_1^2+D_2)\T\,\, \cup\,\, (D_1^2- D_2)\T
\ee
It is an easy matter \cite{ard2} to show that this relation is in fact equivalent to the KP equation
(\ref{hirota}). It is interesting to observe that in (\ref{weil}) only $D_1$ and $D_2$ are involved, while $D_3$ plays no role.
\vskip 0.5 cm
Finally, we come to the  interpretation of the KP equation (\ref{hirota}) in terms of the Kummer map 
\be
\label{kum}
\kappa: X\to |2\T|^*=\mathbb P^N\,,\qquad N=2^g-1\,.
\ee
To say that the image of a point $b\in X$, via the Kummer map is an inflectionary point for the Kummer variety $\kappa(X)$, is like saying that there is a line $l\subset \mathbb P^N$ such that the preimage $\kappa^{-1}(l)$
contains 
the  length 3 artinian subscheme $b+Y\subset X$ associated to some second order germ  $Y$
\be
Y:\,\, \epsilon\mapsto \e  2U+\e^22V\,\,\subset \,
\,X\,.
\ee
We set 
\be
V_Y=\{a=2b\in X\,\,|\,\, b+Y\subset\kappa^{-1}(l)\,\,\text{for some line}\,\,l\subset \mathbb P^N  \}\,\,\subset \,\,X\,.
\ee
Clearly $V_Y\supset Y$ and
for a general abelian variety $X$ the subscheme $V_Y$ could simply coincide with $Y$. On the other hand,
 using Riemann's bilinear relations for second order theta-functions, one can see \cite{ard2} that
 the KP equation (\ref{hirota}) is equivalent to the statement that $V_Y$ contains a \it third order germ \rm extending $Y$:
 \be
 \label{flex}
 V_Y\supset Z:\,\, \epsilon\mapsto \e  2U+\e^22V+\e^32W\,\,\subset \,
\,X\,.
 \ee
On the other hand, if $X=J(C)$ is a Jacobian then $V_Y$ is nothing but the Abel-Jacobi image of $C$ in $J(C)$: much more than a tiny germ
\cite{gun}, \cite{wel1}.
As Novikov conjectured and Shiota proved [S]  the KP equation characterizes Jacobians among principally polarized abelian varieties. This theorem can be stated in the following way.

\vskip 0.2 cm
\sl Let $(X,\T)$ be an indecomposable, principally polarized abelian variety. Then $X$ is the Jacobian of a curve of genus $g$
if and only if there exist vectors $U$, $V$, $W$ in $\mathbb C^g$  (or equivalently constant vector fields $D_1$, $D_2$, $D_3$ on $X$) such that
one of the following equivalent conditions holds.

\begin{enumerate}
\item[(i)\,\,] The KP equation (\ref{kp1}) is satisfied with $u$ as in (\ref{utheta}),
\item[(ii)\,] The system (\ref{eigen}) is satisfied with $u$ as in (\ref{utheta}),
 and $\psi$ as in (\ref{bak}),
\item[(iii)] The KP equation in Hirota's form (\ref{hirota}) is satisfied,
\item[(iv)\,] The Weil relation (\ref{weil}) is satisfied,
\item[(v)\,\,] The Kummer variety of $X$ admits a second order germ of inflectionary tangent(i.e. (\ref{flex})
is satisfied).
\end{enumerate}
\rm Again, we notice that in  (v) the vector $W$ (or equivalently the vector field $D_3$)
makes no appearance.
\vskip 0.5 cm
To state the result of the present paper we go back to the system (\ref{eigen}) and we consider \it only \rm  the first of the two equations:
\be
\label{lequ}
\left(\frac{\partial^2}{\partial x^2}-\frac{\partial}{\partial y}+u\right)\psi=0
\ee
Given $a\in X\setminus\{0\}$ we  look for solutions of (\ref{lequ}) given by
\be
\label{upsia}\,\,\,\,\,\,u=2\frac{\partial^2}{\partial x^2}\log\te(xU+yV+z)\,,\quad
\psi=e^L\cdot
\frac{\te(xU+yV+a+z)}{\te(xU+yV+z)}
\ee
Where $L=Ax+By$. We next express equation (\ref{lequ}) in terms of $\te$ and we get   the bilinear equation
\be
\label{p.AB} \aligned D_1^2\te\cdot\te_a+  \te\cdot D_1^2\te_a+
D_2\te\cdot\te_a-\te\cdot D_2\te_a-2D_1\te\cdot D_1\te_a \\
+2AD_1\te_a\cdot\te-2A\te_a\cdot D_1\te+(A^2-B)\te\cdot\te_a=0\,,
\endaligned
\ee
where $\te_a(z)=\te(z+a)$. Changing $D_2$ into $D_2+AD_1$, we get
\be
\label{p} D_1^2\te\cdot\te_a+  \te\cdot D_1^2\te_a+
D_2\te\cdot\te_a-\te\cdot D_2\te_a-2D_1\te\cdot D_1\te_a +c\te\cdot\te_a=0\,.
\ee
This equation looks much simpler than (\ref{hirota}).
Using the methods we mentioned above, it is  straightforward to show that this equation is equivalent to either of the following Weil-type
relations  \cite{m1}
\be
\label{weil1}
\T\,\,\cap\,\,\T_a\,\,\subset\,\,D_1\T\,\,\cup\,\,D_1\T_a
\ee
\be
\label{weil2}
\qquad\qquad D_1\T\,\,\subset\,\,(D_1^2+D_2)\T\,\,\cup\,\,\T_a
\ee
From the point of view of flexes of the Kummer variety the equation (\ref{p}) simply says that
the there exists a point $b\in X\setminus\{0\}$, with $2b=a$,  such that $\kappa(b)$ is a flex for the Kummer variety $\kappa(X)$, or equivalently that
\be
\label{wa}
a=2b\in V_Y
\ee
It is natural to ask if these equivalent conditions are sufficient to characterize Jacobians among
all principally polarized abelian varieties. This question, in its formulation (\ref{wa}), is a particular case of the so-called \it trisecant conjecture, \rm first formulated in \cite{wel1} (see also \cite{deb}).

In the present paper we give an affirmative answer to this question under the additional hypothesis that
the closure  $<a>$ of the group generated by $a$ is irreducible.

\begin{theorem}
Let $(X,\T)$ be an indecomposable, principally polarized abelian variety. Then $X$ is the Jacobian of a curve of genus $g$
if and only if there exist vectors $U\neq0$, $V$ in $\mathbb C^g$  (or equivalently constant vector fields $D_1\neq0$, $D_2$, on $X$) and a point $a\in X\setminus\{0\}$, with $<a>$ irreducible, 
such that
one of the following equivalent conditions holds.

\begin{enumerate}
\item[(a)] The  equation (\ref{lequ}) is satisfied with $u$ and $\psi$ as in (\ref{upsia}),
\item[(b)] The equation (\ref{p}) is satisfied,
\item[(c)] Either one of  Weil-type relation (\ref{weil1}),(\ref{weil1}) is satisfied,
\item[(d)] There is a point  $b\in X\setminus\{0\}$, with $2b=a$, such that  $\kappa(b)$ is a  flex of the Kummer variety $\kappa(X)$ (i.e. (\ref{wa})
is satisfied).
\end{enumerate}
\end{theorem}
Observe that the ''only if'' part of this theorem is clear:  suppose $X=J(C)$ is the Jacobian of a curve $C$,
and take a general point $b\in\frac{1}{2}\varphi(C)$. Then, on the one hand, the image point $\kappa(b)$ is a flex of the Kummer variety
$\kappa(X)$ and, on the other, $<2b>=X$. In  \cite{kri3} the result, in its formulation \sl (a), \rm is proved under a different hypothesis: namely, that the vector $U$ {\it spans an elliptic curve}.

\vskip 1cm \section {$D_1$-invariant flows.}\label{S: $D_1$-invariant flows.}

In this section we will explain a dichotomy that was first proved in   \cite{m1}.
The dichotomy is the following.
\begin{dichotomy}Let $(X,\T)$ be an indecomposable principally polarized abelian variety.
Assume that equation (\ref{p}) holds for some $a\in X\setminus\{0\}$. Then either the KP equation (\ref{hirota}) holds or the subscheme
$D_1\T\subset\T$ contains a $D_1$-invariant component.
\end{dichotomy}
In the setting of the KP equation, a similar dichotomy was implicit in the work of Shiota \cite{shi} and was also observed in \cite{ard2}.
In that setting the dichotomy was: assume that the KP equation  (\ref{hirota}) holds. Then either the entire KP hierarchy is satisfied or  the subscheme
$D_1\T\subset\T$ contains a $D_1$-invariant component.
\vskip 0.3 cm
Let us  prove  the above dichotomy. From now on we will write 
\be
P=D_1^2\te\cdot\te_a+  \te\cdot D_1^2\te_a+
D_2\te\cdot\te_a-\te\cdot D_2\te_a-2D_1\te\cdot D_1\te_a +c\te\cdot\te_a
\ee
so that equation (\ref{p})  reads:
\be
\label{pp}P=0\,.
\ee
Assume that equation (\ref{pp}) holds and that the  KP equation (\ref{hirota}) does not. In view of the equivalence of (iii) and (iv) in the previous section we may assume that there is an irreducible component $W$ of the subscheme $D_1\T\subset \T\subset X$, such that
\be
\label{notkp}(D_1^2+D_2)\te\cdot(D_1^2-D_2)\te_{|_{W}}\neq 0
\ee
Let $p$ be a general point of the reduced scheme $W_{\text{red}}$.  A theorem of Ein-Lazarsfeld \cite{el} asserts that the theta divisor of an indecomposable abelian variety is smooth in codimension 1.
Since $W_{\text{red}}$ is a divisor in $\T$ the point  $p$ is a smooth point of $\T$. Hence there exist an irreducible element $h$, invertible elements $\beta$, $\gamma$,  elements $\tilde\alpha$ $\tilde\beta$, $\tilde\gamma$,  in the local ring 
$\mathcal O_{X,p}$ and integers $m\geq1$, $r$, $s$ such that the ideal of $W_{\text{red}}$ at $p$ is $(h,\te)$, and such that
\be
\label{locexpr}
\aligned
D_1\te&=h^m+\tilde\alpha\te\,,\\
\te_a&=\beta h^s+\tilde\beta\te\,,\\
D_2\te&=\gamma h^r+\tilde\gamma\te\,.\\
\endaligned
\ee
In particular the ideal of $W$ is $I(W)=(h^m, \te)$. Our goal is to prove that $h$ divides $D_1(h)$.
In fact, in this case, $h$ divides $D_1^n(h)$ for every $n$, so that the $D_1$-line through $p$ is contained in $W\subset D_1\T$. We proceed by contradiction and \it we assume that $h$ does not divide $D_1(h)$. \rm

The next remark is that either $r=0$ or $r\geq m$. For this, use (\ref{locexpr}) to write the identity
$D_1D_2\te=D_2D_1\te$, $\mod (h^r,\te)$. We get $r\gamma D_1(h)h^{r-1}-mD_2(h)h^{m-1}=0$, modulo $(h^m,\te)$.
Since $\gamma$ is invertible, $m$ can not exceed $r$, unless $r=0$. Now write (\ref{notkp}), modulo $(h^m,\te)$, in terms of the local expressions (\ref{locexpr}). As $\gamma$ is invertible, we get that, if $m\geq 2$, then $2r<m$. 
In conclusion, either $r=0$ or $r\geq m=1$.

Look at the equation (\ref{lequ}), where $u$ and $\psi$
are as in (\ref{upsia}). Consider a general point $\eta U+yV+z$ in the theta divisor $\T$.
We  have expansions
\be
\label{psiex}
\psi=\frac{\alpha}{x-\eta}+\beta+\gamma(x-\eta)+\delta(x-\eta)^2+\dots\,,
\ee

\be
\label{uiex}
-\frac{2}{(x-\eta)^2}+v+w(x-\eta)+\dots\,.
\ee
We look at $\eta$, $\alpha$, $\beta$, $\gamma$, $\delta$, $v$, $w$, as function of $y$ and we may assume $\alpha\neq 0$.
Write equation (\ref{lequ})  looking at the coefficients of $(x-\eta)^i$, for $i=-2, -1, 0$. We get
\be
\aligned
\alpha\overset\cdot\eta+2\beta=0\,,\\
-\overset\cdot\alpha+\alpha v-2\gamma=0\,,\\
-\beta+\gamma\overset\cdot\eta+\beta v +\alpha w=0\,.
\endaligned
\ee
Taking the derivative of the first equation and using the last two equations, we get 
\be
\overset{\cdot\cdot}\eta=-2w\,.
\ee
We compute $w$ by recalling the expression of $u$ given in (\ref{upsia}), and $\overset{\cdot\cdot}\eta$ 
by using the identity $\te(\eta(y)U+yV+z)\equiv 0$. We then  obtain the equation
\be
\label{longeq}
\aligned
-D_1^2\te\cdot(D_2\te)^2+2D_1D_2\te\cdot D_2\te\cdot D_1\te-D_2^2\te\cdot(D_1\te)^2=\\
-(D_1^2\te)^3+2D_1^2\te\cdot D_1^3\te\cdot D_1\te-D_1^4\te\cdot(D_1\te)^2\,,\endaligned
\ee
\it which  is valid on $\T$. \rm We plug in (\ref{longeq})   the local  expressions given in (\ref{locexpr}). As we already noticed,  either $r=0$, or $r\geq m=1$. In the first case we look at (\ref{longeq}), modulo $I(W)=(h^m,\te)$ 
and we get 
\be
\left(\gamma^2-m^2D_1(h)^2h^{2m-2}\right)mh^{m-1}D_1(h)=0\,,\qquad\mod (h^m,\te)\,.
\ee
Since $\gamma$ is invertible we must have $m=1$. On the other hand,
The hypothesis (\ref{notkp}) tells us that
$\gamma^2-D_1(h)^2\neq0$, modulo $(h^m,\te)$. It follows that 
$h$ divides $D_1(h)$, proving the dichotomy in this case. If  $r\geq m=1$, we look at (\ref{longeq}), modulo $(h,\te)$ and we get $D_1(h)^3=0$, modulo $(h,\te)$, which again implies that $h$ divides $D_1(h)$.
\hfill Q.E.D.

\vskip 1cm 
\section {The proof of the Theorem}\label{S: Proof of the Theorem}

We keep the notation of the preceding section. We assume that equation (\ref{p}) holds for the theta function of $X$. We consider the divisor $D_1\T\subset \T\subset X$. We say that an irreducible  component $W$ of $D_1\T$ is \it bad \rm,
if (\ref{notkp}) is satisfied.  
We consider the local expressions (\ref{locexpr}).
\begin{Lemma} Let $X$ be an indecomposable principally polarized abelian variety.
Let $a\in X\setminus\{0\}$.
Assume that (\ref{p}) is satisfied.
Then an irreducible component $W$ of $D_1\T$ is bad if and only if 
\be
\label{bad}
W_{\text{red}}\quad \text{is}\quad D_1\,\,\text{- invariant}\,,\quad\text{and}\quad 2r<m
\ee
Moreover,  if  $W$ is bad then $s>m$.
\end{Lemma}
\sl Proof. \rm Suppose that $W$ is bad. By the Dichotomy we know that $ W_{\text{red}}$ is $  D_1 $-invariant or, equivalently, that $D_1 (h) \in (h,\, \theta)$.
Therefore $ D_1^n \theta \in (h^m,\, \theta) \, , \ \forall n \ge 0$.
From (\ref{notkp}) we get: $ 2 r < m$.
The opposite implication is trivial. Now suppose that  $D_1 (h) \in (h,\, \theta)$ and $2r<m$. Computing $ P(\theta) $ modulo $ (h^{\min\{2m, \, m+s\}},\, \theta)$
the only term that survives is $D_2 \theta \cdot \theta_a$.  On the other hand, this product is in
$ (h^{r+s}, \, \theta) $ and therefore $r+s \ge \min\{2m, \, m+s\}$.
Since $2r < m$, we  obtain $r+s \ge 2m$.  In particular $s > m$. \hfill Q.E.D.

\begin{Lemma} Let $X$ be an indecomposable principally polarized abelian variety.
Let $a\in X\setminus\{0\}$. Assume that $<a>$ is irreducible and that
 (\ref{p}) is satisfied. Let $W$ be a bad component of $D_1\T$. Then
$W$ is $a$-invariant.
\end{Lemma}
\sl Proof. \rm
We are going to use the following standard notation.
Let
$ Y$  be a reduced, irreducible variety which is  non-singular in codimension 1. Let
$ Z $ be  an irreducible divisor in $ Y $ and let $h $ be a generator for the ideal of 
$ Z_{\text{red}}$ at a general point  $p\in Z_{\text{red}}$.  
Let $f$ be a regular function at $p$ then we define the symbol  $  o [f ; \, Z ; \, Y ]$
by the identity of ideals in the local ring $ \mathcal O_{p, \, Y}  $
\be
(f) =  \left(  h^{o [f ; \, Z ; \, Y ]}  \right)\,.
\ee
Let now $ \mathcal W$ be the set of bad  irreducible components of $D_1\T$. We claim that
if $W$ belongs to $ \mathcal W$, then $ W_{-a}$ belongs to $ \mathcal W$ as well.
By (\ref{locexpr}) $ \,\,\,m =  o [D_1 \theta; \, W ; \, \Theta ]$,  $\,\,\,s  =  o [\theta_a; \, W ; \, \Theta ]\,\,\,$ and 
$\,\,\,r  = \ o [D_2\theta; \, W ; \, \Theta ] $.
We have $  0 \leq 2 r <  m < s$  and $  o [D_1 \theta_a; \, W ; \, \Theta ]=m  < s$.  As a consequence,
\be
o [D_1 \theta; \, W_{-a} ; \, \Theta ]  =  o [D_1 \theta_a; \, W ; \, \Theta_a ]  =  o [D_1 \theta_a; \, W ; \, \Theta ]  =  m
\ee
and
\be
o [D_2 \theta; \, W_{-a} ; \, \Theta ]  =  o [D_2 \theta_a; \, W ; \, \Theta_a ]  =  o [D_2 \theta_a; \, W ; \, \Theta ]  =  r
\ee
In view of Lemma 1, this proves that $  W_{-a}  \in \mathcal W$. Consider now the irreducible abelian variety $<a>$.
As a consequence  we get: $ \bigcup_{n \in \mathbb N} \, W_{-na} \, \subseteq \, D_1 \Theta $.  Thus,
taking the closure, we have:
\be
\bigcup_{x \in <a>} \, W_x \, \subseteq \, D_1 \Theta \,.
\ee
For dimensional reasons, we conclude that $ W $ is $<a>$-invariant. \hfill Q.E.D.

\sl Proof of Theorem 1. \rm we will finish the proof of Theorem 1 by showing that if $X$ is indecomposable and if its theta function satisfies equation (\ref{p}) with $<a>$ irreducible then $D_1\T$ has no bad component. 
The basic result we   need is the following Lemma, which is reminiscent of Shiota's Lemma B in \cite{shi}.
\begin{Lemma}
Let $ (S, \mathcal L) $ be a polarized abelian variety.
Let $ Y $ be a 2-dimensional disk with analytic coordinates $t $ and $ \lambda $ and let
$\theta $ be a non-zero section of $ \mathcal O_Y \otimes H^0( S, \mathcal L)$.
Let $a$  be a point of $S\setminus \{0\}$. Assume that $<a>$ is irreducible. Let $ D_1 \neq 0$, $\tilde D_2 \in  T_0(S)$
and assume that $ S=\langle D_1 , \, a \rangle $.
Assume that
\be
\label{lb.1}
\aligned
P \theta \quad = \quad
& D_1^2 \theta \cdot \theta_a - 2 D_1 \theta \cdot D_1 \theta_a + \theta \cdot D_1^2 \theta_a + \\
& + \big(\partial_t + \tilde D_2\big) \theta \cdot \theta_a - \theta \cdot \big(\partial_t + \tilde D_2\big) \theta_a +
c \cdot \theta \cdot \theta_a \quad = \quad 0 \ ,
\endaligned
\ee

\be
\label{lb.2}
\theta (t, \lambda, x) \quad = \quad
\sum_{i, j \geq 0} \theta_{i,j}(x) \cdot t^i \lambda^j \ , \qquad x\in S\,.
\ee
Write:
\be
\theta (t, \lambda, x) \quad = \quad t^\nu\lambda^\rho\left[\theta_{\nu,\rho}(x)+t\alpha(t, x)\right]+\lambda^{\rho+1}\beta(t,\lambda,x)\,,
\ee
where $ \theta_{\nu,\rho}\not\equiv 0$.
Furthermore, assume $ \nu \ge 1  $. Then there exist local sections at zero 
$f\in \mathcal O_Y $ and $ \psi\in\mathcal O_Y \otimes H^0(  S, \mathcal L )$ such that
\be
\label{dec.t.lambda}
\theta (t, \lambda, x) \quad = \quad \lambda^{\rho} \cdot f(t,\lambda) \cdot \psi(t,\lambda,x) \ ,
\ee
where $ \psi (0,0,\cdot ) \not\equiv 0 , \ f(0,0) = 0 \ $ and $  f( \cdot , 0 ) \not\equiv 0$.
\end{Lemma}
\vskip 0.3 cm
It is important to observe that the geometrical meaning of (\ref{dec.t.lambda}) is the the following:
\be
\label{dec.t.lambda.geom}\{\te(t, \lambda, x)=0\}\,\cap\,\,
\{D_1\te(t, \lambda, x)=0\}\,\,\supset\,\,\{\lambda^\rho\cdot f(t, \lambda)=0\}\,.
\ee
\vskip 0.3 cm
We assume this Lemma for the time being and we continue the proof of Theorem 1. 
We proceed by contradiction. We then suppose that a bad component $W$ of $D_1\T$  exists.
Let $S$ be the Zariski closure of the subgroup generated by the $D_1$-flow and the point $ a $.
We shall write $ S  = \langle D_1 , \, a \rangle $ and 
we set $ \mathcal L = \mathcal O(\T)\vert_S $.
Since $ D_1 \neq 0 \ $ we have $S \neq 0 $
on the other hand, by Lemma 2,  $ W $ contains a translate of $ S $ therefore $  S\neq X $.
Note that $ W_{\text{red}} $ is $ \ T_0(S) $-invariant.
Let $B$ be the complement of $S$ in $X$ relative to the polarization $\Theta $.
In the sequel we shall work on $B \times S$
via the natural isogeny $\pi : B \times S \rightarrow X$.
We shall also write $\theta$ instead of $ \pi^{\star} \theta$ 
while working on $B \times S$. 
Our abuse of notation reflect the fact that $ \theta $ and
$ \pi^{\star} \theta$ coincide as theta functions
on the common universal cover of $X$ and $B \times S$. 
Let us fix a general point $p$ of $W_{\text{red}}$.
Clearly, up  a translation of $\Theta$,  we are free to assume $p = 0\in
X$.
Let $\mathcal R= ( W \cap B )_{\text{red}}$.
Observe that $ W_{\text{red}}$ is the $ T_0(S)$-span of $\mathcal R$, that
is:
  $W_{\text{red}}= \mathcal R+ S$. Also observe  that $\mathcal R$ has
codimension 2 in $B$.

Let us decompose $D_2$ as $\tilde D_2 + \overline D_2$,
where $\tilde D_2 \in T_0 ( S )$ and $ \overline  D_2 \in T_0 (B) $.
Since $t_x^\star \theta\in \, (h^m,\, \theta)\,, \forall x \in S$, then
$ \tilde D_2 \theta \in (h^m,\, \theta)$. On the other hand $D_2
\theta \not\in (h^m,\, \theta)$.
It follows that $\overline D_2 \neq \, 0$.
Let $L $ be the analytic germ at zero of the $\overline  D_2
$-integral line in $B$ through zero,
let $\mathcal C$ be the germ at $0$ of a smooth curve in
$B $ meeting $L$ transversally only at $ 0$, and let $Y$ be the
surface $\mathcal C+ L $ in $B$.
Let $ \lambda $ be a parameter on $ \mathcal C$ vanishing at $0$ and
let $ t$ be the coordinate on $L$,
vanishing at zero, with $ \partial_t = \overline D_2 $.
Thus $\lambda$ and $t $ are parameters on $Y$. Set $ \Omega = Y
\times  S\subset  B \times S$.
On $\Omega $ we write
\be
\label{theta.t.l}
\theta (t,\lambda,x) =
\sum_{i, j \geq 0} \ t^i \cdot \lambda^j \cdot \theta_{i,j}(x) \ ,
\qquad x  \in  S \,.
\ee
As $B$ and $S$ are complementary with respect to $ \Theta $,
the $\theta (t, \lambda, \cdot ) $'s, as well as their derivatives
$\theta_{i,j} = {1 \over i! \cdot j! } \left( { \partial^j \over
\partial \lambda^j } { \partial^i \over
\partial t^i } \theta \right) (0,0,\cdot) \ $,
belong to $ \ H^0( S, \Theta \vert_{S} ) $.

 Our analysis will distinguish  two cases  corresponding to
whether the variety $ \mathcal R$ is  $ \overline D_2 $-invariant, or not.

\vskip 0.5 cm
Let us first assume that $ \mathcal R$ is not $ \overline D_2 $-invariant.
In this case, to reach a contradiction, our gaol will be twofold. On the one hand,
we will choose $\mathcal C$ in such a way that
\be
\label{l=t=0}
\Omega \cap \pi^{-1}(W_\text{red})=  \{ \lambda = t = 0 \} \times S\,.
\ee
On the other, we will show that $\te(t,\lambda,x)$ satisfies the hypotheses of Lemma 3
with $\rho=0$:
$$
\theta(t,\lambda,x) = f(t, \lambda) \cdot \psi ( t, \lambda, x ) \ ,
$$
where $f(0,0) = 0 $.  The conclusion will be that
\be
\Omega \cap \pi^{-1} W  = \Omega \cap \{ \theta = 0 \} \cap \{ D_1 \theta = 0 \} \, \supseteq \, \Omega \cap \{ f = 0 \}\,.
\ee
But since $f(0,0) = 0$,  it would  follow that $ \Omega \cap \pi^{-1} W $ has codimension 1 in
$ \Omega $ contradicting (\ref{l=t=0}).

To achieve our goals,
we choose $ \mathcal C$ in such a way that it meets $ \mathcal R$  transversally only at
$0$ and that $\partial_{\lambda}$ does not belong to  $\langle T_0 (\mathcal R) , \overline D_2 \rangle \, \cup \, T_0(\Theta)$.
This is possible because $\mathcal R$ has codimension 2 in $B$.
Note that the locus $\{ t = 0 \}$ which is $ \mathcal C\times S $,  is transverse to $ \Theta$.
As a consequence, we are free to assume that the function $t $ is the restriction $h \vert_{Y \times S}$.
We have $Y \cap \mathcal R = \{ \lambda = t = 0 \}$.  Thus (\ref{l=t=0}) holds.
Observe now, that the two restrictions $\Theta \vert_{Y\times S}$ and $\{ t = 0 \} \vert_\Omega$,
meet along $\Omega \cap  \overline  W_\text{red} $ and they are smooth at the general point of
$\Omega \cap  \overline  W_\text{red} $. It follows that $ \theta_{i,0}\vert_S \not \equiv 0 $ for some $ i$. On the other hand,
since $\T\supset\{0\}\times S$, we must have
$ \theta_{0,0}\vert_S = 0 $.
 We can therefore apply lemma 3 with $\rho=0$.

\vskip 0.5 cm

Let us now assume that $ \mathcal R$ is $ \overline D_2 $-invariant. In this case we will choose
$  \mathcal C $ in such a way that 
\be
\label{l=0}
 \Omega \cap \pi^{-1} ( W_{\text{red} }) = \{ \lambda = 0 \}\,.
 \ee
At the same time we will prove that, also in this case, Lemma 3 applies so that:
$$
\theta (t,\lambda,x) \quad = \quad \lambda^{\rho} \cdot f(t, \lambda) \cdot \psi ( t, \lambda, x ) \, , \
$$
so that $ f $ divides both $  \theta\vert_{\Omega} \ $ and $ D_1 \theta\vert_{\Omega}$.
It would then follow that
\be
\label{case2}
\Omega \cap \pi^{-1} ( W_{\text{red}} ) \, \supset \, \Omega \cap \{ f = 0 \}\,.
\ee
We would also have$  f(0,0) = 0 $ and $ f(\cdot,0)  \not\equiv 0 $.
But then (\ref{case2}) would tell us that the locus $  \Omega \cap \pi^{-1} ( W_{\text{red}} )  $ contains, locally at $ p$,  a component which
is not the component $  \{ \lambda = 0 \}$
contradicting (\ref{l=0}).

To follow this line of reasoning, we
choose $  \mathcal C $ in such a way that it is contained in $\Theta $ and meets $ \mathcal R$
transversally only at $0$.
Since the loci $ \{ h = 0 \} $ and $ \Theta $ are transverse and $ { \mathcal C } $ meets $\mathcal R$ transversally at $ 0$,
we may assume that $  \lambda  $ is the restriction of $ h $ to $  \mathcal C  \times \{ s \} \cong { \mathcal C }$. Thus (\ref{l=0}) holds.
Now write
\be
\label{case2}
\theta(t,\lambda,x)  = \lambda^\rho \cdot \sum_i \, t^i \cdot \theta_{i,\rho}(x) \ + \ \sum_{i, \, j > \rho} \ t^i \cdot \lambda^j \cdot \theta_{i,j}(x)\,, 
\ee
The hypotheses of Lemma 3 require us to show  that $
\theta_{0,\rho}(x) \equiv0$.
Since $S$ is generated by $a$ and by the flow of $D_1$,
it suffices to prove that $D_1^i\theta_{0,\rho}(na)=0$ for all $i$
and $n$ in $\mathbb N$.
Now, on one hand we have
\be
D_1^i\theta_{na}(0,\lambda, 0) \ = \ \lambda^\rho \cdot
D_1^i\theta_{0,\rho}(na) \ , \quad \mod \lambda^{\rho+1} \, .
\ee
On the other hand, with the same notation from (\ref{locexpr}) and lemma 1,
we can write
$ D_1^i\theta_{na} = \delta h^m + \tilde\delta \theta $,
furthermore,
since $ D_2 $ does not involve $\lambda$ we have that $ r \ge \rho$,
so that, again in view of lemma 1, we have $m>\rho$.
Working modulo $ \lambda^{\rho+1} $ we get
$ D_1^i\theta_{na}(0,\lambda,0) =
(\delta h^m + \tilde\delta \theta)\vert_{t=x=0} =
\tilde\delta(0,\lambda,0)\theta(0,\lambda,0) = 0 \, , \ $
where the last equality holds since $ \mathcal C$ is
contained in $\Theta $.
This proves that $D_1^i\theta_{0,\rho}(na)=0$ as required.

\vskip 0.5 cm
This ends the proof of Theorem 1. It remains to prove Lemma 3.

\vskip 0.5 cm
\sl Proof of Lemma 3. \rm
Set
$ \omega(t,x) =  \sum_{i \geq \nu} \theta_{i,\rho}(x) \cdot t^{i-\nu}\,\,$  so that
\be
\aligned
 \theta &= \lambda^\rho t^\nu \cdot \omega\,,\qquad \mod ( \lambda^{\rho+1} )\,,\\
 \omega (0,x)  &= \theta_{\nu,\rho} (x) \quad \not \equiv \quad 0 \,.
 \endaligned
\ee
The reader is advised to follow the computations by setting $\nu=1$. The notation needed for the general case,
somewhat overwhelms the reasoning.
We first construct the function $f(t, \lambda)$ and the section $\psi(t, \lambda, x)$ as formal power series
in $t$ and $\lambda$. To this end,
we look for constants and sections

\be
c_{i,j}\in \mathbb C \ , \quad 0 \leq i \leq \nu-1 \, , \ j \ge 1 \ , \quad
g_{i,j} (x)  \in  H^0 (S , \mathcal L ) \ , \quad i \geq \nu \, , \ j \ge 1 \ ,
\ee
such that
\be\label{lb.3}
\theta (t,\lambda,x) =  \lambda^\rho \cdot \left( t^\nu \ + \ \sum_{j\geq 1 } f_j (t) \cdot \lambda^j \right)
\cdot \left( \omega(t,x)  +  \sum_{j \geq 1} \omega_j (t,x) \cdot \lambda^j \right) \,,
\ee
where, for $  j \geq 1$,  we define
\be
\label{lb.4}
f_j (t)  =  \sum_{i=0}^{\nu-1} c_{i,j} \cdot t^i \ , \qquad \omega_j (t,x)  =  \sum_{i\geq \nu} g_{i,j}(x) \cdot t^{i-\nu} \,.
\ee
Define  $  \tilde P (r,s)  =  \frac{ 1}{2} [ P(r+s) - P(r) - P(s) ]$.
It is  straightforward to verify the following properties:

\be 
\label{lb.5.1}  P (r)  =  \tilde P (r,r)\,,\quad\text{and} \quad\tilde P\,,\quad\text{is a symmetric} \quad\Bbb C [\lambda] -\text{bilinear operator}\,,
\ee

\be
\label{lb.5.2}  \tilde P ( g \cdot r , g \cdot s )  =  g^2 \cdot \tilde P (r ,s )\,,\quad \text{for}\quad g = g(t,\lambda)\,\quad
\text{not depending on}\quad x\,,
\ee

\be
\label{lb.5.3}\tilde P ( t^i \cdot r , t^j \cdot s ) = 
t^{i+j} \cdot \tilde P (r,s) + {1 \over 2}(i-j)t^{i+j-1} \cdot (r \cdot s_u - r_u \cdot s) \, . 
\ee
In particular, $  P ( \lambda^{\rho} \cdot r )  =  \lambda^{ 2 \rho} \cdot P (r)  $ and therefore we are free to assume
$\rho = 0$.  Furthermore, writing $  \theta  =  t^\nu \cdot \omega + \lambda  \theta'  $ one has

\be
0  =  P(\theta)=  t^{2\nu} \tilde P\left(\omega ,  \omega \right)\,,\quad \mod (\lambda)\,.
\ee
As a consequence,
\be
\label{lb.6}
P \left( \omega(t,x) \right)= 0\,.
\ee
We now proceed by induction.  Let $ k $ be a positive integer, and assume that we found
constants $  c_{i,j}$,  for all $  1 \leq j \leq k-1$, $i \leq \nu-1$,
and sections $  g_{i,j}(x)$,  for all $  1 \leq j \leq k-1$, $i \geq \nu$,
such that (\ref{lb.3}) holds modulo $ (\lambda^k)$.  Set
\be
g(t,\lambda)  =  t^\nu + \sum_{j=1}^{k-1} f_j(t) \cdot \lambda^j \, , \qquad
\phi(t,\lambda,x)  =  \omega + \sum_{j=1}^{k-1} \omega_j(t,x) \cdot \lambda^j\,,
\ee
and define $  \omega'(t,x) $ by
\be
\label{lb.7}
\theta =  g \cdot \phi + \lambda^k \cdot \omega' \, , \qquad \mod  (\lambda^{k+1}) \, .
\ee
We need to prove that there exist constants $  c_{i,k} \, , \ i \leq \nu-1$,
and sections $  g_{i,k}\, , \,\, i \geq \nu$, such that
\be
\label{lb.8}
\omega' (t,x)  = 
\sum_{i=0}^{\nu-1} c_{i,k} \cdot t^i \cdot \omega(t,x) 	+  \sum_{i \geq \nu } g_{i,k}(x) \cdot t^i \, .
\ee
In fact, defining $ f_k ,  \omega_k  $ as (\ref{lb.4}) requires, it is clear that (\ref{lb.3}) holds modulo $ ( \lambda^{k+1} )$.
Working modulo $  (\lambda^{k+1})  $ one has
\be
0  = P(\theta)= \ g^2 \cdot P(\phi) + 2 \lambda^k \tilde P( t^{\nu} \cdot \omega , \omega') \,,\quad\mod (\lambda^{k+1}) \,.
\ee
In particular we get $  g^2 \cdot P(\phi) = 0$,  modulo  $(\lambda^k)$.
Since $  g^2(t,\lambda)  =  t^{2\nu}$, modulo  $(\lambda)$,   is non-zero, we get $  P(\phi) = 0$, modulo $(\lambda^k)$.
Since $  g^2 \cdot P(\phi) + 2 \lambda^k \tilde P( t^{\nu} \cdot \omega , \omega')  =  0$, modulo $(\lambda^{k+1})$,
and (again) $  g^2(t,\lambda)  =  t^{2\nu}$ modulo $(\lambda)  $, we get
\be
\label{lb.9}
\tilde P ( t^{\nu} \cdot \omega, \omega') =  0\,, \qquad \mod (t^{2\nu}) \, .
\ee
To prove (\ref{lb.8}) we now proceed by induction on $ i $.  We assume there exists $ i_0$,  satisfying
$ 0 \leq i_0 \leq \nu-1$,   such that
\be
\label{lb.10}
\omega' (t,x)  =  \sum_{i=0}^{i_0-1} c_{i,k} \cdot t^i \cdot \omega (t,x)  +  \eta (x) \cdot t^{i_0} \, ,\quad
\mod ( t^{i_0 + 1} ) \, ,
\ee
and we have to prove that $ \eta(x) $ is a multiple of $  \omega (t,x) $ modulo $(t)$.
Equivalently, we have to prove that $ \eta(x) $ is a multiple of $  \omega (0,x)$.
Since $  \tilde P (\omega,\omega) = P ( \omega ) = 0 $, (see (\ref{lb.6})), applying (\ref{lb.5.3}) we get
$  \tilde P ( t^{\nu} \cdot \omega , t^i \cdot \omega )  =  0$.
Now, substituting (\ref{lb.10}) in (\ref{lb.9}) and using again (\ref{lb.5.3}), we get the following equality  modulo $(t^{\nu+i_0})$:
\be
\aligned
0  =& \tilde P \left( t^{\nu} \cdot \omega \, , \ \sum_{i=0}^{i_0-1} c_{i,k} \cdot \omega \cdot t^i + \eta (x) \cdot t^{i_0} \right)\\
= & (\nu-i_0) \cdot t^{\nu+i_0-1} \cdot \omega (0,x) \cdot \omega (0,x-a) \cdot
\left( { \eta (x) \over \omega (0,x) } - { \eta (x-a) \over \omega (0,x-a) }\right) \, . 
\endaligned
\ee
It follows that the meromorphic function on $S$
\be
 c_{i_0,k}(x) = \frac{\eta (x)}{\omega (0,x) } 
 \ee
 is invariant under  translation by $ a  $.
 We want to show that $c_{i_0,k}(x)$   is constant.
We proceed by contradiction. Suppose it is not.
As $ c_{i_0,k}(x) $ is $ a $-invariant, its poles are $ a $-invariant, so that
the zero locus of $ \omega (0,x) $ contains 
an $a$-invariant divisor $U$. Since $<a>$ is irreducible we may well assume that $U$ is irreducible.
We want to show that
\be
\label{lemma.A}
D_1^\alpha \omega \vert_{U}  =  0 \, , \quad \forall  \alpha \in \mathbb N \, .
\ee
We assume (\ref{lemma.A}) and we postpone for the moment its  proof. As $  S $ is generated by  $a$ and the $D_1$ flow,
(\ref{lemma.A}) gives $\omega(0,x)\equiv0$.
This is a contradiction.
Thus, $  \eta (x) =  c_{j_0,k} \cdot \omega (0,x)$, so that
$ \omega' (t,x)  =  \sum_{i=0}^{i_0} c_{i,k} \cdot t^i \cdot \omega (t,x)$, modulo  $( t^{i_0 + 1} ) $.

At this stage both $f$ and $\psi$ are constructed as formal power series. Now we prove that
they are in fact regular.
As
$  \psi (0,0,\cdot) \not\equiv 0  $ we are allowed to fix a point  $  x_0  $ such that
$  \psi (0,0,x_0) \neq 0  $ and consider the formal power series
$  q(t, \lambda)  $ defined by 
$  \psi (t,\lambda,x_0) \cdot q(t, \lambda)  =  1$.
Set
$  \tilde f (t, \lambda) =  f (t, \lambda) \cdot \psi
(t,\lambda,x_0)  $
and $  \tilde \psi (t,\lambda,x) =  \psi (t,\lambda,x) \cdot q(t, \lambda)$.
It is clear that
$  \theta (t, \lambda, x)  =  \lambda^{\rho} \cdot \tilde f (t, \lambda)
\cdot \tilde \psi (t, \lambda, x)$.
As $ \tilde \psi (t,\lambda,x_0)  =  1  $ and
$  \theta (t,\lambda,x_0)  $ are both convergent,
$ \tilde f (t, \lambda) \ $ is  convergent as well.
But now the convergence of $  \tilde \psi (t,\lambda,x)  $ follows from the one of $ \theta (t,\lambda,x) $ and 
$ \tilde f (t, \lambda) $.
Note that $  t^{\nu} \ $ divides  $  \tilde f (t,0) = 0$ , that
$\tilde f (t,0) \not \equiv 0  $
and that $  \tilde \psi (0,0,\cdot) \not \equiv 0  $.
Thus the properties we need hold for 
 $ \tilde f$ and $ \tilde \psi $.
 \vskip 0.5 cm
To finish the proof of Lemma 3, it now remains to prove (\ref{lemma.A}). By hypothesis $\omega(t,x)\in\mathcal O_\Delta\times H^0(\mathcal L,S)$, where $\Delta$ is a 1-dimensional disc centered at the origin and with coordinate $t$. Also $\omega$
 satisfies 
 \be
\label{pomega}
\aligned
P \omega \quad = \quad
& D_1^2 \omega \cdot \omega_a - 2 D_1 \omega \cdot D_1 \omega_a + \omega \cdot D_1^2 \omega_a + \\
& + D_2 \omega \cdot \omega_a - \omega \cdot D_2\omega_a +
c \cdot \omega \cdot \omega_a =  0 \ ,
\endaligned
\ee
where $D_2=\partial_t+\tilde D_2$, while  $D_1$ and $\tilde D_2$ are constant vector fields on $S$.
Moreover there is an irreducible $a$-invariant divisor $U$ in $\{0\}\times S$ on which $\omega$ vanishes.
Under these hypotheses we want to prove (\ref{lemma.A}). 
We proceed by contradiction. We assume that there exists $ b > 0 $
such that $  D_1^b \omega \vert_U  \not \equiv \ 0 $.
Let
\be
\label{la.1}
\aligned
\beta_{\gamma} & = \text{ min}  \{  \beta  \,\,\vert\,\,  D_1^{\beta} D_2^{\gamma} \omega \vert_U  \not \equiv  0  \} \, , \\
w  & = \text{ min}  \{ \beta_{\gamma} + 2 \gamma  \} \, ,\\
\sigma  & =  \text{ max}  \{  \gamma \,\, \vert\,\, \beta_{\gamma} + 2\gamma = w  \} \, ,\\
\endaligned
\ee
Observe that $  w  \leq  \beta_0  \leq  b  <  \infty$,
$\,\,  \sigma \leq \frac{1 }{2} w < \infty $,
$ \,\,w = \beta_{\sigma} + 2 \sigma$ and
$\,\, w \leq \beta_{\gamma} + 2 \gamma  $ for all $\,\, \gamma $.
As $\,  \omega \vert_U \equiv 0  $ we have $ \beta_0 > 0 $.
Moreover
$ \ D_1^{\beta} D_2^{\gamma} \omega \vert_U  =  0 , \ $ for all $  \beta < \beta_{\gamma} $.
It follows that
\be
\label{la.2}
\aligned
\text{if}\quad \beta +2 \gamma < w \, ,  &\quad \text{then} \quad D_1^{\beta} D_2^{\gamma} \omega \vert_U  = 0\,, \\
\text{if}\quad \gamma > \sigma \quad \text{and} \quad \beta +2 \gamma \leq w \, , &\quad \text{then} \quad D_1^{\beta} D_2^{\gamma} \omega \vert_U  =  0\,.
\endaligned
\ee
Let now
\be
\label{la.3}
\aligned
\tilde w  & =  \text{ min}  \{  \beta_{\gamma} + 2\gamma \,\, \vert \,\, \gamma<\sigma  \} \,,\\
\tilde \sigma  & =  \text{ max}  \{ \gamma  \,\,\vert  \,\,\gamma<\sigma ,  \beta_{\gamma} + 2\gamma = \tilde w  \}\,.
\endaligned
\ee
 Since  $w=\sigma+\beta_\sigma\leq  \tilde w=\tilde\sigma+\beta_{\tilde\sigma}$, we have  $\beta_{\tilde\sigma} \geq 2 $. Then,
\be
\label{la.4}
\aligned
0  & = D_1^{\beta_{\tilde\sigma} - 2} D_2^{\tilde\sigma+\sigma} P \, \omega \, \vert_U \\
& = {{\tilde\sigma + \sigma} \choose {\sigma}}
\left( D_1^{\beta_{\tilde\sigma}} D_2^{\tilde\sigma} \omega \cdot D_2^\sigma \omega_a \, + \,
D_1^{\beta_{\tilde\sigma}} D_2^{\tilde\sigma} \omega \cdot D_2^\sigma \omega_a \right) \vert_U \ .
\endaligned
\ee
Now, observe that $  D_1^{\beta'} D_2^{\gamma'} \omega \vert_U \, \in \, H^0( U, \mathcal L ) $, provided that lower order derivatives
vanish on $ \, U \, , \ $ i.e.: $ D_1^\beta D_2^\gamma \omega \vert_U \, = \, 0  $ for
$  \beta \leq \beta'$, $\gamma \leq \gamma'$,  $\beta + \gamma < \beta' + \gamma'$.
In particular,
\be
D_1^{\beta_{\tilde\sigma}} D_2^{\tilde\sigma} \omega \vert_U \ , \quad D_2^\sigma \omega \vert_U \qquad \in \qquad H^0( U, \mathcal L) \, .
\ee
On the other hand, this two sections are non-trivial, therefore we have on $U$ a non zero meromorphic function
$  f  := \, { D_1^{\beta_{\tilde\sigma}} D_2^{\tilde\sigma} \over D_2^\sigma \omega} \vert_U$ which,
by (\ref{la.4}),  satisfies $ f \, = \, - f_a$.  As $ <a>$  is irreducible and of positive dimension,
there exists a sequence of odd multiples of $ a$ converging to the origin of $ <a>$. As $U$ is $<a>$-invariant,
for general $ x \in U $  there exists a sequence of points $  x_i  $ converging to $  x $  and such that 
$ f(x) = - f(x_i)$.  This gives $ f = 0 $ which is a contradiction.\hfill Q.E.D.

\vskip 1cm 
\section {Final remarks}\label{S: Final remarks}

We end this note by connecting equation (\ref{p}), or better  (\ref{p.AB}), with the KP hierarchy.
As is well known (see for instance formula (3.29) in \cite{ard2}) the KP hierarchy for the theta function can be written
in the form
\be
\label{kp.hierar}
\aligned
 \e\left(D_1^2\te\cdot\te_{2\zeta(\e)}+  \te\cdot D_1^2\te_{2\zeta(\e)}+
D_2\te\cdot\te_{2\zeta(\e)}-\te\cdot D_2\te_{2\zeta(\e)}-2D_1\te\cdot D_1\te_{2\zeta(\e)}\right) \\
+D_1\te_{2\zeta(\e)}\cdot\te-\te_{2\zeta(\e)}\cdot D_1\te+d(\e)\te\cdot\te_{2\zeta(\e)}=0\,,
\endaligned
\ee
where 
\be
\zeta(\e)=\e U+\e^2V+\e^3W+\dots\,,\quad d(\e)=d_3\e^3+d_4\e^4+\dots
\ee
The similarity of this equation with (\ref{p.AB}) is obvious. In fact  write  (\ref{p.AB}) with 
\be
\label{germs}
a=a(\e)=2\zeta(\e)\,, \quad A=A(\e)=\frac{1}{\e}+ \e^2(a_0+a_1\e+\dots)\,,\quad B(\e)=\frac{1}{\e^2}+ \e(b_0+b_1\e+\dots)\,,
\ee
with $2a_0=b_0$. Then change parameter from $\e$ to $\frac{1}{A}$, multiply the resulting equation by this parameter
to get exactly the KP hierarchy (\ref{kp.hierar}).

In view of the equivalence of (\ref{lequ}) and (\ref{p.AB}), the KP hierarchy can also be expressed by saying that equation  (\ref{lequ})holds, where now
\be
\label{upsiL}
u=2\log\frac{\partial^2}{\partial x^2}\te(xU+yV+z)\,,\quad
\psi=e^L\cdot
\frac{\te(xU+yV+a(\e)+z)}{\te(xU+yV+z)}\,,
\ee
and $L=A(\e)x+B(\e)y$ and $a(\e)$, $A(\e)$ and $B(\e)$ as in (\ref{germs}).

\vskip 1 cm

\footnotesize{\sl Acknowledgments. \rm We wish to thank Sam Grushevski for very useful conversation on the subject of this paper. The first named 
 author wishes to thank Columbia University, The Italian Academy, and the Accademia dei Lincei for  hospitality and support during the preparation of this work.}

\vskip 1 cm

\def\cprime{$'$}

\bibliography{Flex}
\bibliographystyle{alpha}
\end{document}